\documentstyle[12pt,leqno]{article}
\begin{document}
\title{Locally Lagrangian Symplectic and Poisson Manifolds}
\author{{\normalsize by}\\Izu Vaisman}
\date{}
\maketitle
{\def\thefootnote{*}\footnotetext[1]%
{{\it 2000 Mathematics Subject Classification}
53D05, 53D12, 53D17. \newline\indent{\it Key words and phrases}:
Tangent Structures, Lagrangian Symplectic
Structures, Lagrangian Poisson Structures}}
\begin{center} \begin{minipage}{12cm}
A{\footnotesize BSTRACT. We discuss symplectic manifolds where, locally, the
structure is that encountered in Lagrangian dynamics. Examples and
characteristic properties are given. Then, we refer to the computation of
the Maslov classes of a Lagrangian submanifold. Finally, we indicate the
generalization of this type of symplectic structures to Poisson manifolds.
The paper is the text of a lecture presented at the Conference ``Poisson
2000" held at CIRM, Luminy, France, between June 26 and June 30, 2000. It
reviews results contained in the author's papers
\cite{{V1},{V2},{V3}} as well as in papers by other authors.}
\end{minipage} \end{center}
\section{Locally Lagrangian Symplectic Manifolds}
The present paper is the text of a lecture presented at the Conference
``Poisson
2000" held at CIRM, Luminy, France, between June 26 and June 30, 2000, and
it reviews results contained in the author's papers
\cite{{V1},{V2},{V3}} as well as in papers by other authors. The notion of
a locally Lagrangian Poisson manifold is defined for the first time here.

The symplectic structures used in Lagrangian dynamics are defined on a
tangent manifold $TN$, and consist of symplectic forms of the type
$$\omega_{\cal L}=\frac{1}{2}\left(\frac{\partial^{2}{\cal L}}
{\partial q^{i}\partial
u^{j}}-\frac{\partial^{2}{\cal L}}{\partial q^{j}\partial
u^{i}}\right)dq^{i}\wedge dq^{j}+\frac{\partial^{2}{\cal L}}{\partial u^{i}
\partial u^{j}}du^{i}\wedge dq^{j}. \leqno{(1.1)}$$
In (1.1), $(q^{i})_{i=1}^{n}$ $(n=dim\,N)$ are local coordinates in the
configuration space $N$, $(u^{i})$ are the corresponding natural
coordinates in the
fibers of $TN$, and ${\cal L}$ is the non degenerate Lagrangian function on
$TN$ (i.e., ${\cal L}\in C^{\infty}(TN)$, $rank(\partial^{2}{\cal L}/
\partial u^{i}\partial u^{j})=n$). (In this paper, everything is
$C^{\infty}$.)

The most known geometric description of $\omega_{\cal L}$ is that it is the
pullback of the canonical symplectic form of $T^{*}N$ by the Legendre
transformation defined by ${\cal L}$. But, geometrically, it is more
significant that $\omega_{\cal L}$ is related with the {\em tangent
structure} of the manifold $TN$. The latter is the bundle morphism
$S:TTN\rightarrow TTN$ defined by
$$SX\in T{\cal
V},\;(SX)(<\alpha_{\pi(u)},u>)=<\alpha_{\pi(u)},\pi_{*}X_{u}>,
\leqno{(1.2)}$$
where ${\cal V}$ is the foliation of $TN$ by fibers (the {\em vertical
foliation}), $u\in TN$,  $X\in\Gamma TTN$,
$\alpha\in\Gamma T^{*}N$, and $\pi:TN\rightarrow N$ is the natural
projection. ($\Gamma$ denotes spaces of global cross sections.)
Formulas (1.2) define the action of $SX$ on $q^{i},u^{i}$, and one has
$$S\left(\frac{\partial}{\partial q^{i}}\right)
=\frac{\partial}{\partial u^{i}},
\;S\left(\frac{\partial}{\partial u^{i}}\right)=0.\leqno{(1.3)}$$
The relation between $\omega_{\cal L}$ and $S$ is e.g., \cite{JK}
$$\omega_{\cal L}=d\theta_{\cal L},\;\theta_{\cal L}=d{\cal L}\circ S.
\leqno{(1.4)}$$

A symplectic form given by (1.4) is called a {\em global Lagrangian
symplectic (g.L.s.) structure}. The Lagrangian ${\cal L}$ on $TN$ uniquely
defines $\omega_{\cal L}$ but, if (and only if) we change ${\cal L}$ by
$${\cal L}'={\cal L}+f(q)+\alpha_{i}(q)u^{i},\leqno{(1.5)}$$
where $f\in C^{\infty}(N)$ and $\alpha_{i}(q)dq^{i}$ is a closed $1$-form
on $N$, we get the same symplectic form $\omega_{{\cal L}'}=\omega_{\cal
L}$.

Accordingly, it is natural to study locally Lagrangian symplectic
manifolds $M$ i.e., manifolds which have an open covering
$M=\cup_{\alpha\in A}U_{\alpha}$, with local Lagrangian functions ${\cal
L}_{\alpha}\in C^{\infty}(U_{\alpha})$ that give rise to a global
symplectic form which locally is
of the type (1.4). In particular, this is interesting since
we may expect to also have compact manifolds of this type while, a g.L.s.
manifold cannot be compact since its symplectic form is exact.

First, we must ask our manifold $M$ to carry a {\em tangent structure}
\cite{CB}
i.e., an {\em almost tangent structure} $S\in\Gamma End(TM)$ where
$$S^{2}=0,\;\;im\,S=ker\,S,\leqno{(1.6)}$$
which is {\em integrable} i.e., locally defined by (1.3). The existence of
$S$ implies the even-dimensionality of $M$, say $dim\,M=2n$, whence
$rank\,S=n$, and integrability is equivalent with the vanishing of the
Nijenhuis tensor
$${\cal N}_{S}=[SX,SY]-S[SX,Y]-S[X,SY]+S^{2}[X,Y]=0.\leqno{(1.7)}$$

It is important to notice that $S$ defines the {\em vertical distribution}
$V:=im\,S=ker\,S$, which is tangent to a foliation ${\cal V}$ if $S$ is
integrable. Moreover, in the latter case ${\cal V}$ is {\em locally
leafwise affine} since a change of local coordinates where (1.3) holds is
of the form
$$\tilde q^{i}=\tilde q^{i}(q^{j}),\;\;\tilde u^{i}=\frac{\partial\tilde
q^{i}}{\partial q^{j}}u^{j}+\theta^{i}(q^{j}).\leqno{(1.8)}$$
(We use the Einstein summation convention.)The parallel vector fields of
the locally affine structure of the leaves of ${\cal V}$ are the vector
fields $SX$ where $X$ is a ${\cal V}$-projectable vector field on $M$.

Now, we can give the formal definition: a {\em locally Lagrangian
symplectic (l.L.s.) manifold} is a manifold $M$ endowed with a tangent
structure $S$ and a symplectic structure $\omega$ of the form (1.4), where
${\cal L}$ are local functions on $M$.

A simple example can be obtained as follows. Take $N={\bf R}^{n}$,
$TN={\bf R}^{2n}$, and the Lagrangian of the {\em modified harmonic
oscillator} \cite{Md}
$${\cal
L}=\frac{1}{2}(\delta_{ij}u^{i}u^{j}+\alpha_{ij}q^{i}q^{j}),\leqno{(1.9)}$$
where $(\delta_{ij})$ is the unit matrix and $(\alpha_{ij})$ is a constant
symmetric matrix. Then, quotientize by the group
$$q^{i}\mapsto q^{i}+m^{i},\;u^{i}\mapsto u^{i}+s^{i}\hspace{3mm}
(m^{i},s^{i}\in{\bf Z})$$
to get the torus $T^{2n}$. The tangent structure of ${\bf R}^{2n}$ projects
to $T^{2n}$, and the function ${\cal L}$ yields local functions on $T^{2n}$
which have transition relations (1.5) hence, define the same Lagrangian
symplectic form. Notice that this example is on a compact manifold. Further
examples will be given later on.

Following is a coordinate-free characterization of the l.L.s. manifolds
\cite{V3}:
\proclaim 1.1 Proposition. Let $M$ be a manifold with a tangent structure
$S$ and a symplectic form $\omega$. Then, $\omega$ is locally Lagrangian
with respect to $S$ iff
$$\omega(X,SY)=\omega(Y,SX)\hspace {5mm}(\forall X,Y\in\Gamma
TM).\leqno{(1.10)}$$ \par \noindent{\bf Proof.}
Using local coordinates where (1.3) and (1.1) hold, it is easy to check
(1.10). Conversely, from (1.10), it follows that the vertical foliation
${\cal V}$ of $S$ ($T{\cal V}=im\,S$) is $\omega$-Lagrangian. Hence, $M$
may be covered by local charts $(U,x^{i},y^{i})$ such that
$$V:=T{\cal V}=span\left\{\frac{\partial}{\partial
y^{i}}\right\},\;\;\omega=\sum_{i}dx^{i}\wedge dy^{i}.\leqno{(1.11)}$$
Furthermore, we must also have
$$S\left(\frac{\partial}{\partial
x^{i}}\right)=\sum_{k}\lambda_{ik}\frac{\partial}{\partial y^{k}},
\leqno{(1.12)}$$
where, because of (1.10), $\lambda_{ik}=\lambda_{ki}$.

Now, (1.7) implies $$[S\frac{\partial}{\partial x^{i}},
S\frac{\partial}{\partial x^{j}}]=0$$
hence, there are new coordinates
$$q^{i}=x^{i},\;u^{i}=u^{i}(x^{j},y^{j}) \leqno{(1.13)}$$
such that
$$S\left(\frac{\partial}{\partial q^{i}}
\right)=S\left(\frac{\partial}{\partial x^{i}}
+\frac{\partial y^{k}}{\partial q^{i}}\frac{\partial}{\partial y^{k}}
\right)=S\left(\frac{\partial}{\partial x^{i}}
\right)=\frac{\partial}{\partial u^{i}}.
\leqno{(1.14)}$$

From (1.12) and (1.14) we get $\partial y^{i}/\partial u^{k}=\lambda_{ki}$,
and, since $\lambda_{ki}=\lambda_{ik}$, there exist local functions ${\cal
L}$ such that $y^{i}=-\partial{\cal L}/\partial u^{i}$. Using (1.11), we
get (1.1). Q.e.d.

(1.10) is the {\em compatibility condition} between $S$ and $\omega$. It
allows us to notice one more interesting object namely,
$$\Theta([X]_{V},[Y]_{V}):=\omega(SX,Y),\leqno{(1.15)}$$
where the arguments are cross sections of the transversal bundle $\nu{\cal
V}:=TM/V$ of the foliation ${\cal V}$. $\Theta$ is a well defined
pseudo-Euclidean metric with the local components $(\partial^2{\cal
L}/\partial u^{i}\partial u^{j})$. If this metric is positive definite, we
say that the manifold $(M,S,\omega)$ is of the {\em elliptic type}.
\proclaim 1.2 Proposition. Let $(M,\omega)$ be a symplectic manifold
endowed with a Lagrangian foliation ${\cal V}$ $(T{\cal V}=V)$, and a
${\cal V}$-projectable pseudo-Euclidean metric $\Theta$ on $\nu{\cal V}=
TM/V$. Then, there exists a unique $\omega$-compatible tangent structure
$S$ on $M$ for which $\Theta$ is the metric (1.15).\par
\noindent{\bf Proof.} Split $TM=V'\oplus V$, where $V'$ also is
a $\omega$-Lagrangian distribution, and identify $\nu{\cal V}$ with $V'$.
Then define $$S/_{V}=0,\;\;S/_{V'}=\sharp_{\omega}\circ\flat_{\Theta},
\leqno{(1.16)}$$
where the musical isomorphisms are defined as in Riemannian geometry.
It is easy to check that $S$ is the required tangent structure. In
particular, we must check that ${\cal N}_{S}(X,Y)=0$, and it suffices to
look at the various cases where the arguments are in $V,V'$ while, if in
$V'$, they are ${\cal V}$-projectable vector fields. The only non trivial
case $X,Y\in\Gamma V'$ is settled by noticing that, whenever $X,Y,Z\in
\Gamma_{proj}V'$, one has
$$d\omega(Z,SX,SY)=\omega(Z,[SX,SY])=0,$$
whence $[SX,SY]=0$. Q.e.d.

Proposition 1.2 allows us to find more examples of l.L.s. manifolds.
Let
$$H(1,p):=\{\left(\begin{array}{ccc}Id_{p}&X&Z\\
0&1&y\\0&0&1 \end{array}\right)\;/\;X,Z\in{\bf R}^{p},y\in{\bf R}\}
\leqno{(1.17)}$$
be the {\em generalized Heisenberg group}, and take the quotients
$$M(p,q):=\Gamma(p,q)\backslash(H(1,p)\times H(1,q)),$$
where $\Gamma(p,q)$ consists of pairs of matrices of type (1.17) with
integer entries. Then, with the notation of (1.17) on the two factors, the
form
$$\omega=\,^{t}dX_{1}\wedge(dZ_{1}-X_{1}dy_{1})+
\,^{t}dX_{2}\wedge(dZ_{2}-X_{2}dy_{2})+dy_{1}\wedge dy_{2}\leqno{(1.18)}$$
($t$ denotes matrix transposition) defines a symplectic structure on
$M(p,q)$, and the equations
$$X_{1}=const.,\;X_{2}=const.,\;y_{1}-\alpha y_{2}=const.\hspace {3mm}
(\alpha\in{\bf R}),\leqno{(1.19)}$$
define an $\omega$-Lagrangian foliation ${\cal V}$ \cite{{CFG},{V4}}.
Furthermore
$$g:=\,^{t}dX_{1}\otimes dX_{1}+\,^{t}dX_{2}\otimes dX_{2}+
[d(y_{1}-\alpha y_{2}]\otimes[d(y_{1}-\alpha y_{2}]\leqno{(1.20)}$$
$$+\,^{t}(dZ_{1}-X_{1}dy_{1})\otimes(dZ_{1}-X_{1}dy_{1})
+\,^{t}(dZ_{2}-X_{2}dy_{2})\otimes(dZ_{2}-X_{2}dy_{2})+dy_{2}\otimes dy_{2}
$$
is a projectable metric on the transversal bundle of ${\cal V}$.

Correspondingly, the construction of Proposition 1.2 yields a l.L.s.
structure of elliptic type on $M(p,q)$.

A similar construction holds on the so-called Iwasawa manifolds
$$I(p)=\Gamma_{c}(1,p)\backslash H_{c}(1,p),$$
where $H_{c}(1,p)$ is given by (1.17) with ${\bf R}$ replaced by ${\bf C}$,
and $\Gamma_{c}(1,p)$ consists of matrices with Gauss integers entries
\cite{{CFG},{V4}}.

We end this section by
\proclaim 1.3 Proposition. The l.L.s. manifold $(M,S,\omega)$ is globally
Lagrange symplectic iff $\omega=d\epsilon$ for some global $1$-form
$\epsilon$ on $M$ such that (i) $\epsilon$ vanishes on the vertical leaves
of $S$, and (ii) if $\eta$ is the cross section of $V^{*}$ ($V=T{\cal V}$,
${\cal V}$ is the vertical foliation of $S$) which satisfies $\eta\circ
S=\epsilon$, then $\eta=d_{\cal V}{\cal L}$, where $d_{\cal V}$ is the
differential along the leaves of ${\cal V}$ and ${\cal L}\in
C^{\infty}(M)$. \par \noindent{\bf Proof.}
Notice the exact sequence
$$0\rightarrow V\stackrel{\subseteq}{\rightarrow}TM
\stackrel{\pi}{\rightarrow}TM/V\approx V\rightarrow 0$$
which shows that $V^{*}$ can be identified with the subbundle $(TM/V)^{*}$
of $T^{*}M$. Accordingly, we see that a unique leafwise form $\eta$
as required is
associated with each $\epsilon$ that satisfies (i)
namely, using the local coordinates of
(1.3), we have
$$\epsilon=\epsilon_{i}dq^{i}\mapsto
\eta=\epsilon_{i}[du^{i}]_{(TM/V)^{*}}.
\leqno{(1.21)}$$

If we are in the g.L.s. case, we have (1.4), and $\epsilon=\theta_{\cal L}$
is the required form.

Conversely, if $\omega=d\epsilon$, where (i), (ii) hold, we must locally
have
$$\epsilon=(\frac{\partial{\cal L}}{\partial u^{i}}du^{i})\circ S=
\frac{\partial{\cal L}}{\partial u^{i}}dq^{i}, \leqno{(1.22)}$$
and ${\cal L}$ is the required global Lagrangian. Q.e.d.

Notice that condition (ii) has a cohomological meaning. If
the l.L.s. form $\omega=d\epsilon$ where $\epsilon$
satisfies (i) then (1.4) implies
$$\epsilon/_{U}=\theta_{{\cal
L}_{U}}+\lambda_{U},\hspace{5mm}d\lambda_{U}=0,$$
where $U$ is an open neighborhood on which the local Lagrangian ${\cal
L}_{U}$ exists, and $\lambda$ is a $1$-form which vanishes on the leaves of
${\cal V}$. Thus,
$$\lambda=d(\varphi(q)),\hspace{3mm}
\eta=\left(\frac{\partial{\cal L}_{U}}{\partial u^{i}}+
\frac{\partial\varphi(q)}{\partial q^{i}}\right)[du^{i}],$$ and
we get $d_{\cal V}\eta=0$. In foliation theory (e.g., \cite{V5}), it is
known that if $H^{1}(M,\Phi_{\cal V})=0$, where $\Phi_{\cal V}$ is the sheaf
of germs of functions $f(q^{i})$, then $d_{\cal V}\eta=0$ implies
$\eta=d_{\cal V}{\cal L}$ for a global function ${\cal L}$.
\section{Further geometric results}
a). Let $(M,S)$ be a manifold endowed with an integrable tangent structure
$S$. Then, the following problem is of an obvious interest: find the
symplectic structures $\omega$ on $M$ which are compatible with $S$ hence,
are l.L.s. forms.
\proclaim 2.1 Proposition. The symplectic form $\omega$ on $M$ is
$S$-compatible iff
$$\omega=\frac{1}{2}\varphi_{ij}(q)dq^{i}\wedge
dq^{j}+d(\zeta_{i}(q,u)dq^{i}),\leqno{(2.1)}$$
where $(q^{i},u^{i})$ are the local coordinates of (1.3), the first term of
(2.1) is a local closed $2$-form and $\zeta=\zeta_{i}dq^{i}$ is a local
$1$-form such that $d_{\cal
V}(\zeta_{i}[du^{i}])=0$ and $d\zeta$ is non degenerate. \par
\noindent{\bf Proof.} If $\omega$ is $S$-compatible,
$\omega$ is of the form (1.4), which is (2.1)
with a vanishing first term.

For the converse result, we notice that our hypotheses indeed imply that
the form $\omega$ is closed and non degenerate, while the compatibility
condition (1.10) is equivalent to $\partial\zeta_{i}/\partial u^{j}=
\partial\zeta_{j}/\partial u^{i}$ i.e., $d_{\cal
V}(\zeta_{i}[du^{i}])=0$. Q.e.d.

The local Lagrangians of the form $\omega$ of (2.1) are obtained by putting
(locally)
$$\frac{1}{2}\varphi_{ij}(q)dq^{i}\wedge
dq^{j}=d(\alpha_{i}(q)dq^{i}),\;\;\zeta_{i}[du^{i}]=d_{\cal V}f,\;\;f\in
C^{\infty}(M).\leqno{(2.2)}$$
Then the Lagrangians are
$${\cal L}=f+\alpha_{i}u^{i}.\leqno{(2.3)}$$

A coordinate-free criterion which ensures (2.1) is given by
\proclaim 2.2 Proposition. The symplectic form $\omega$ on $(M,S)$ has the
local form (2.1), possibly without $d_{\cal
V}(\zeta_{i}[du^{i}])=0$, iff the vertical foliation ${\cal V}$ of $S$ is a
Lagrangian foliation with respect to $\omega$. \par
\noindent{\bf Proof.} Consider the vertical foliation ${\cal V}$ of $S$,
denote $V=T{\cal V}$, and let $V'$ be a transversal distribution. We will
use a well known technique of foliation theory namely, the bigrading of
differential forms and multivector fields associated with the decomposition
$TM=V'\oplus V$. In particular, one has
$$d=d'_{(1,0)}+d''_{(0,1)}+\partial_{(2,-1)}, \leqno{(2.4)}$$
where $d''$ may be identified with $d_{\cal V}$, and $d^{2}=0$ becomes
$$d''^{2}=0,\;\partial^{2}=0,\;d'd''+d''d'=0,\;d'\partial+\partial d'=0,
\;d'^{2}+d''\partial+\partial d''=0. \leqno{(2.5)}$$
Furthermore, $d''$ is the coboundary of the leafwise de Rham cohomology,
and it satisfies a Poincar\'e lemma e.g., \cite{V5}.

Clearly, ${\cal V}$ is Lagrangian with respect to (2.1).

Conversely, if the foliation ${\cal V}$ is Lagrangian for a symplectic form
$\omega$, we must have a decomposition
$$\omega=\omega_{(2,0)}+\omega_{(1,1)},\leqno{(2.6)}$$
and $d\omega=0$ means
$$d''\omega_{(1,1)}=0,\;d''\omega_{(2,0)}+d'\omega_{(1,1)}=0,\;
d'\omega_{(2,0)}+\partial\omega_{(1,1)}=0.\leqno{(2.7)}$$

Accordingly, there exists a local $(1,0)$-form $\zeta$ such that
$\omega_{(1,1)}=d''\zeta$ (the $d''$-Poincar\'e lemma), and (with (2.5)) the last
two conditions (2.7) become
$$\begin{array}{ll} d''\omega_{(2,0)}+d'd''\zeta=&d''(\omega_{(2,0)}-
d'\zeta)=0,\vspace{1mm}\\  d'\omega_{(2,0)}+\partial d''\zeta=&
d'\omega_{(2,0)}-d'^{2}\zeta-d''\partial\zeta=d'(\omega_{(2,0)}-d'\zeta)=0,
\end{array} \leqno{(2.8)}$$
whence, $$\Phi:=\omega_{(2,0)}-d'\zeta
\leqno{(2.9)}$$ is a closed $2$-form of bidegree $(2,0)$.
Therefore, $\omega$ has the local expression (2.1). Q.e.d.

The simplest geometric case is that of a tangent bundle $M=TN$ with the
canonical tangent structure (1.3). In this case, the local expressions (2.1)
can be glued up by means of a partition of unity on $N$, and we get all the
global $S$-compatible symplectic forms \cite{V2}
$$\omega=\pi^{*}\Phi+d\zeta,\leqno{(2.10)}$$
where $\pi:TN\rightarrow N$, $\Phi$ is a closed $2$-form on $N$, and
$\zeta$ is a $1$-form on $TN$ which satisfies the hypotheses of Proposition
2.1. In particular, $d_{\cal V}(\zeta_{i}[du^{i}])=0$,
and the contractibility of
the fibers of $TN$ allows to conclude that $\zeta_{i}[du^{i}]=d_{\cal
V}\varphi$, $\varphi\in C^{\infty}(TN)$. This also yields another
expression of the $S$-compatible forms on $TN$:
$$\omega=\pi^{*}\Phi+\omega_{\varphi},\leqno{(2.11)}$$
where $\omega_{\varphi}$ is given by (1.1).

Furthermore, if we also use Proposition 1.3, we see that $\omega$ of (2.10)
is globally Lagrangian iff it is exact. Indeed,if $\omega$ of (2.10)
is $\omega=d\xi$, where $\xi=\xi_{(1,0)}+\xi_{(0,1)}$, we must have
$d''\xi_{(0,1)}=0$ since $\omega$ has no $(0,2)$-component. Then, because
the fibers of $TN$ are contractible, $\xi_{(0,1)}=d''\varphi$, $\varphi\in
C^{\infty}(TN)$, and
$$\omega=d\xi_{(1,0)}+dd''\varphi=d\xi_{(1,0)}+d(d-d')\varphi=
d(\xi_{(1,0)}-d'\varphi)=d\epsilon,$$
where $\epsilon=\xi_{(1,0)}-d'\varphi$ is a $(1,0)$-form. If (2.10) is
l.L.s., then it is clear that $\epsilon$ satisfies the hypotheses of
Proposition 1.3, with a global Lagrangian of the form (2.3).

In this context it is interesting to mention that a criterion to
distinguish the tangent bundles in the class of the manifolds $(M,S)$ with
an integrable tangent structure is available. Namely, $M$ is a tangent
bundle iff the following three conditions hold: (i) the vertical foliation
${\cal V}$ of $S$ has simply connected leaves, (ii) the flat affine
connection induced by $S$ on the leaves of ${\cal V}$ is complete, (iii)
${\cal E}({\cal V})=0$, where ${\cal E}({\cal V})$ is the $1$-dimensional
cohomology class with coefficients in the sheaf of germs of ${\cal
V}$-projectable cross sections of $T{\cal V}$ produced by the difference
cocycle of the local {\em Euler vector fields} $E:=u^{i}(\partial/\partial
u^{i})$ \cite{{Cr},{V4}}.
\vspace{2mm}\\
b). At this point we restrict ourselves to the case of a tangent
bundle $M=TN$ with the canonical tangent structure $S$ of (1.3). Then, we
may speak of {\em second order vector fields} $X$ on $TN$ i.e., vector
fields $X$ such that the projection on $N$ of their trajectories  satisfy an
autonomous system of ordinary differential equations of the second order.
With respect to the local coordinates (1.3) such a vector field has the
form $$X=u^{i}\frac{\partial}{\partial q^{i}}+\gamma^{i}(q,u)
\frac{\partial}{\partial u^{i}}. \leqno{(2.12)}$$

For $TN$, the local coordinate transformations (1.8) of $u^{i}$ are linear
(i.e., $\theta^{i}=0$ in (1.8)) and $E=u^{i}(\partial/\partial u^{i})$ is a
well defined global vector field on $TN$
(the infinitesimal generator of
the homotheties)  called the {\em Euler vector
field}. We see that the vector field $X$ is of the second order iff $SX=E$.
\proclaim 2.3 Proposition. The symplectic form $\omega$ on $TN$ is
$S$-compatible iff the following conditions are satisfied: (i) the vertical
foliation ${\cal V}$ by the fibers of $TN$ is $\omega$-Lagrangian, (ii)
there exist $\omega$-Hamiltonian vector fields which are second order
vector fields. \par
\noindent{\bf Proof.} \cite{V2}. If $\omega$ is $S$-compatible, it must be
of the form (2.11), and (i) holds. Concerning (ii), it is known in
Lagrangian dynamics that, if we consider the energy associated with the
Lagrangian $\varphi$ of (2.11) given by
$${\cal E}_{\varphi}=E\varphi-\varphi,$$
its $\omega_{\varphi}$-Hamiltonian vector field $X_{\varphi}$
($i(X_{\varphi})\omega_{\varphi}=-d{\cal E}_{\varphi}$) is of the second
order.

Furthermore, all the
second order vector fields are given by $X=X_{\varphi}+Z$, where $Z$
is an arbitrary vertical vector field, and the $1$-form
$\Psi=i(X)\pi^{*}\Phi$ ($\Phi$ of (2.11)) is independent on $Z$.

If we want
a function $h\in C^{\infty}(TN)$ such that
$$-dh=i(X)\omega=i(X_{\varphi}+Z)\omega
=\Psi-d{\cal E}_{\varphi}+i(Z)\omega_{\varphi},$$
it means we want a relation of the form
$$\Psi+i(Z)d''\zeta=df,\leqno{(2.13)}$$
where $\zeta$ is the $1$-form of (2.10), $d''$ is defined using an
arbitrary complementary distribution $V'$ of $V=T{\cal V}$, and
$f={\cal E}_{\varphi}-h$ must be the lift of a function on $N$, since the
left hand side of (2.13) is of bidegree $(1,0)$.

Because the Lagrangian $\varphi$ is non degenerate, (2.13) has a solution
$Z$ for any $f\in C^{\infty}(N)$. Therefore, condition (ii) of Proposition
2.3 is satisfied and we even know how to find all the $\omega$-Hamiltonian,
second order vector fields.

Conversely, condition (i) implies (2.10) (see Proposition 2.2), and, if a
vector field $X$ of the form (2.12) such that $i(X)\omega=-dh$ ($h\in
C^{\infty}(TN)$) exists, the equality of the $(0,1)$-components yields
$$u^{i}\frac{\partial\zeta_{i}}{\partial u^{k}}=
\frac{\partial h}{\partial u^{k}} \leqno{(2.14)}$$
whence the derivatives
$$\frac{\partial\zeta_{j}}{\partial u^{k}}=\frac{\partial^{2}h}{\partial
u^{j}\partial u^{k}}-
u^{i}\frac{\partial^{2}\zeta_{i}}{\partial u^{j}\partial u^{k}}$$
are symmetric, and we are done (see Proposition 2.1). Q.e.d.

Concerning second order vector fields on a tangent bundle $TN$, the
following problem is important:
\proclaim Problem 1. Let $X$ be a second order vector field on TN.
Study the existence and generality of the Poisson structures
$P$ on $TN$ such that
$X$ is a $P$-Hamiltonian vector field. \par
Conversely, we can formulate:
\proclaim Problem 2. If $P$ is a Poisson structure on $TN$, study the
existence and generality of second order $P$-Hamiltonian vector fields.
\par
In analogy with the variational calculus problems, Problem 2 can be called
the {\em direct problem}, and Problem 1 the {\em inverse problem} for
Hamiltonian second order vector fields. Proposition 2.3 gives the solution
of the direct problem for symplectic structures $\omega$ on $TN$ for which
the vertical foliation ${\cal V}$ is Lagrangian
\cite{V2}. But, it is easy to
check that if $W$ is a Poisson structure on $N$ its lift $P$ to $TN$
has no second order Hamiltonian vector fields. $P$ is defined by
$\{f\circ\pi,g\circ\pi\}=0$ for all $f,g\in C^{\infty}(N)$, $\pi:
TN\rightarrow N$, and
$$\{\alpha,f\circ\pi\}=(\sharp_{W}\alpha)f,\;
\{\alpha,\beta\}=L_{\sharp_{W}\alpha}\beta-L_{\sharp_{W}\beta}\alpha-
d(W(\alpha,\beta)),$$
for fiberwise linear functions identified with $1$-forms $\alpha,\beta$.
We again refer to \cite{V2} for an iterative method of solving the inverse
problem for Poisson structures $P$ such that $\{f\circ\pi,g\circ\pi\}=0$.

It is interesting to notice that the notion of a second order vector field
may also be defined for arbitrary manifolds endowed with an integrable
almost tangent structure.

Let $(M,S)$ be a manifold with the integrable structure (1.3).
Then, (1.8) shows the existence of a natural foliated structure on the
vector bundle $V=T{\cal V}$ (${\cal V}$ is the vertical foliation of $S$),
and it makes sense to speak of ${\cal V}$-projectable cross sections
of $V$. A vector
field $X\in\Gamma TM$ will be called a {\em second order vector field} if
for each canonical coordinate neighborhood $(U,q^{i},u^{i})$, $SX-E_{U}$,
where $E_{U}=u^{i}(\partial/\partial u^{i})$ is the local Euler vector
field, is a ${\cal V}$-projectable cross section of $V$.
The condition is invariant by the coordinate
transformations (1.8), and the local expression of $X$ is of the form
$$X=(u^{i}+\alpha^{i}(q))\frac{\partial}{\partial q^{i}}
+\beta^{i}(q,u)\frac{\partial}{\partial u^{i}},\leqno{(2.15)}$$
which yields
$$\frac{d^{2}q^{i}}{dt^{2}}=\beta^{i}(q^{j},\frac{dq^{j}}{dt}
-\alpha^{i}(q))+\frac{\partial\alpha^{i}}{\partial q^{j}}
\frac{dq^{j}}{dt} \leqno{(2.16)}$$
along the trajectories of $X$. Equation (2.16) explains the name.

Notice that, if ${\cal V}$ is a fibration, second order vector fields
always exist. It suffices to glue up local second order fields by a
partition of unity on the basis.

The following fact, which is known for tangent bundles (e.g. \cite{Md})
is also true in the general case: if $X$ is a second order vector field then
$(L_{X}S)^{2}=Id$, where $L$ denotes the Lie derivative. This is easily
checked by using (2.15) and by acting on $\partial/\partial q^{i},\,
\partial/\partial u^{i}$. It follows that $F:=L_{X}S$ is an almost product
structure on $M$ which has the $(+1)$-eigenspace equal to $V=T{\cal V}$.
The complementary distribution
$$V'=span\left\{\frac{\partial}{\partial q^{i}}-\frac{1}{2}
\left(\frac{\partial\alpha^{j}}{\partial q^{i}}
-\frac{\partial\beta^{j}}{\partial u^{i}}\right)\frac{\partial}{\partial
u^{j}}\right\} \leqno{(2.17)}$$
is the $(-1)$-eigenspace of $F$.

In particular, if the functions $\beta^{i}$ of (2.15) are quadratic with
respect to $u^{i}$ (a condition which is invariant by (1.8)), $V'$ is an
{\em affine distribution} transversal to ${\cal V}$ i.e., the process of
lifting paths of transversal submanifolds of ${\cal V}$ to paths tangent to
$V'$ yields affine mappings between the leaves of ${\cal V}$
\cite{{Mol},{V4}}.
\vspace{2mm}\\
c). In what follows, we give a result concerning symplectic
reduction of locally Lagrangian symplectic manifolds.
\proclaim 2.4 Proposition. Let $(M,S,\omega)$ be a l.L.s. manifold and $N$
a coisotropic submanifold with the kernel foliation
$C=(TN)^{\perp_{\omega}}$. Assume that the following conditions are
satisfied: $$S(TN)\subseteq TN,\;\;V\cap TN\subseteq S(TN)+C\hspace{3mm}
(V=im\,S),\leqno{(i)}$$
(ii) the leaves of $C$ are the fibers of a submersion $r:N\rightarrow P$,\\
(iii) the restriction of $S$ to $N$ sends $C$-projectable vector fields to
$C$-projectable vector fields.\\
Then $S$ induces a tangent structure $S'$ on $P$ which is compatible with
the reduction $\omega'$ of the symplectic structure $\omega$ to $P$.\par
\noindent{\bf Proof.} (i) and the compatibility condition (1.10) imply
$S(C)\subseteq C$, therefore, we get an induced morphism
$\tilde S:TN/C\rightarrow TN/C$, such that $\tilde S^{2}=0$ and
$im\,\tilde S=(V\cap TN)/(V\cap C)$.

Since the quotient is Lagrangian for the reduction of $\omega$,
we get $$rank\,\tilde S=\frac{1}{2}dim\,(TN/C).$$

Condition (ii) allows us to reduce $\omega$ to a symplectic structure
$\omega'$ of $P$, and
(iii) ensures that $\tilde S$ projects to an $\omega'$-compatible
tangent structure
$S'$. Q.e.d.
\section{Maslov Classes}
Since a l.L.s. manifold $(M,S,\omega)$ has a canonical Lagrangian foliation
${\cal V}=im\,S$, any Lagrangian submanifold $L$ of $M$ has Maslov classes,
which are cohomological obstructions to the transversality of $L$ and
${\cal V}$. In this section, we will indicate a differential geometric way
of computation of these Maslov classes.

We begin by a brief recall of the general
definition of the Maslov classes \cite{V1}. Let $\pi:E\rightarrow M$ be a
vector bundle of rank $2n$, and $\sigma\in\Gamma\wedge^{2}E^{*}$ a nowhere
degenerate cross section. Then $(E,\sigma)$ is a {\it
symplectic vector bundle}. Furthermore, $(E,\sigma)$ has a reduction
of its structure group $Sp(2n,{\bf R})$ to the unitary group $U(n)$
which is defined up to homotopy, and can be fixed by the choice of a
complex structures $J$ {\em calibrated} by $\sigma$
(i.e., $\sigma$-compatible:
$\sigma(Js_{1},Js_{2})=\sigma(s_{1},s_{2})$, and such that
$g(s_{1},s_{2}):=
\sigma(s_{1},Js_{2})$ is positive definite, $s_{1},s_{2}\in\Gamma E$).

If $L_{0}\rightarrow M$ is a Lagrangian subbundle of $(E,\sigma)$,
the complex version of $(E,J,g)$ has
local unitary bases of the form
$$\epsilon_{i}=\frac{1}{\sqrt{2}}(e_{i}-\sqrt{-1}Je_{i}), \leqno{(3.1)}$$
where $(e_{i})_{i=1}^{n}$ is a local, real, $g$-orthonormal basis of $L_{0}$
($L_{0}$-{\em orthonormal}, $J$-{\em unitary bases}).

Correspondingly, $E$ admits $L_{0}$-orthogonal, $J$-unitary connections
$\nabla^{0}$ of
local expressions
$$\nabla^{0}\epsilon_{i}=\theta_{0i}^{j}\epsilon_{j}\hspace{1cm}
(\theta_{0i}^{j}+\theta_{0j}^{i}=0), \leqno{(3.2)}$$
and a local curvature matrix
$$\Theta^{0}=d\theta^{0}+\theta^{0}\wedge\theta^{0}.\leqno{(3.3)}$$

Furthermore, if $L_{1}\rightarrow M$ is a second Lagrangian subbundle
of $(E,\sigma)$, and
$\nabla^{1}$ is a $L_{1}$-orthogonal, $J$-unitary connection,
the following objects exist:\\
i) the {\it difference tensor}
$\alpha:=\nabla^{1}-\nabla^{0}$,\\
ii) the {\it curvature variation}
$$\Theta_{t}:=(1-t)\Theta_{0}+t\Theta_{1}+t(1-t)\alpha\wedge\alpha,
\leqno{(3.4)}$$
$0\leq t\leq1$ (of course, $\alpha$ and $\Theta$ are matrices),\\
iii) the {\it Chern-Weil-Bott forms}
$$\Delta(\nabla^{0},\nabla^{1})c_{2h-1}:=\leqno{(3.5)}$$
$$=(-1)^{h+1}\frac{\sqrt{-1}}{(2\pi)^{2h-1}(2h-2)!}
\int_{0}^{1}(\delta_{i_{1}...i_{2h-1}}^{j_{1}...j_{2h-1}}\tilde\alpha^{i_{1}}
_{j_{1}}\wedge\tilde\Theta_{tj_{2}}^{i_{2}}\wedge...\wedge
\tilde\Theta_{tj_{2h-1}}^{i_{2h-1}})dt, $$
where the components
$\tilde\alpha^{i}
_{j}$, $\tilde\Theta_{tj}^{i}$ of $\alpha$ and $\Theta_{t}$ are taken with
respect to any {\it common}, local, $J$-unitary bases of $(E,J,g)$
(which may not be $L_{a}$-orthogonal $(a=0,1)$).

It turns out that the forms (3.5) are closed, and define cohomology classes
$$\mu_{h}(E,L_{0},L_{1}):=[\Delta(\nabla^{0},\nabla^{1})c_{2h-1}]
\in H^{4h-3}(M,{\bf R}), \leqno{(3.6)}$$
($h=1,2,...$), which do not depend on the choice of $J$, $\nabla^{0}$,
$\nabla^{1}$, and are invariant by homotopy deformations of $(L_{0},L_{1})$
via Lagrangian subbundles.

We refer to \cite{V1} for details. The classes (3.6) are called the
{\it Maslov classes} of the pair $(L_{0},L_{1})$. If $L_{0}\oplus L_{1}=E$,
one may use $\nabla^{1}=\nabla^{0}$ and $\mu_{h}=0$.
The following two properties are also important:
$$\mu_{h}(E,L_{0},L_{1})=-\mu_{h}(E,L_{1},L_{0}), \leqno{(3.7)}$$
$$\mu_{h}(E,L_{0},L_{1})+\mu_{h}(E,L_{1},L_{2})+
\mu_{h}(E,L_{2},L_{0})=0.\leqno{(3.8)}$$

If $(M,\omega)$ is a symplectic manifold endowed with a Lagrangian
foliation ${\cal V}$, and if $L$ is a Lagrangian submanifold of $M$, we
have {\em Maslov classes} of $L$ defined by
$$\mu_{h}(L)=\mu_{h}(TM/_{L},T{\cal V}/_{L},TL))\in H^{4h-3}(L,{\bf R}).
\leqno{(3.9)}$$
For $h=1$, and if $L$ is a Lagrangian submanifold of ${\bf
R}^{2n}=T^{*}{\bf R}^{n}$, $\mu_{1}(L)$ is the double of the original
class defined by Maslov.

In particular, if $(M,S,\omega)$ is a l.L.s. manifold, and ${\cal V}$ is
its vertical Lagrangian foliation, formula (3.9) will be the
definition of the Maslov classes of the Lagrangian submanifold $L$ of $M$.
We will discuss a way of computing these classes.

Let $V'$ be a transversal Lagrangian distribution of $V=im\,S$. Then
$F:=S/_{V'}$ is an isomorphism $V'\approx V$ with the inverse
$F^{-1}:V\approx V'$, and it is easy to check that, if $F^{-1}$ is extended
to a morphism $S':TM\rightarrow TM$ by asking $S'/_{V'}=0$, then $S'$ is an
almost tangent structure on $M$ which satisfies the compatibility condition
$$\omega(X,S'Y)=\omega(Y,S'X).\leqno{(3.10)}$$
Furthermore, it is also easy to check that $J:=S'-S$ is an
$\omega$-compatible almost complex structure, which is positive iff
$(M,S,\omega)$ is of the elliptic type (i.e., $\omega(SX,Y))$ is a positive
definite bilinear form on $V'$; see Section 1).

Therefore, on l.L.s. manifolds of the elliptic type there is an easy
construction of an almost complex structure $J$ as needed in the
computation of the Maslov classes, and we also get the corresponding
Riemannian metric
$$g(X,Y):=\omega (X,JY).\leqno{(3.11)}$$
We will say that $g$ is the {\em Riemannian metric associated with} $V'$,
and the restriction of $g$ to $V'$ is the metric $\Theta$ defined in
Section 1.

Now, we need connections as required in (3.5). We can obtain such
connections in the following way. Start with a metric connection
$\nabla^{0}$ of the vector bundle $V$, and extend $\nabla^{0}$ to $TM$ by
asking
$$\nabla^{0}_{X}Z=S'\nabla^{0}_{X}Y,\leqno{(3.12)}$$
if $Z=S'Y\in\Gamma V'$ ($Y\in\Gamma V$).
It follows easily that the extended connection, also denoted by
$\nabla^{0}$, is a $V$-orthogonal, $J$-unitary connection, and we will use
it on $TM/_{L}$.

Furthermore, since $L$ is a Lagrangian submanifold and $J$ is compatible
with $\omega$, the $g$-normal bundle of $L$ is $JTL$, and we may write down
{\em Gauss-Weingarten equations} of the form \cite{V1}
$$\begin{array}{lcr} \nabla^{0}e_{i}&=&\lambda^{j}_{i}e_{j}+b^{j}_{i}
(Je_{j}),\vspace{1mm}\\ \nabla^{0}(Je_{i})&=&-b^{j}_{i}e_{j}+
\lambda_{i}^{j}(Je_{j}).\end{array} \leqno{(3.13)}$$
In (3.13), $(e_{i})$ is a $g$-orthonormal basis tangent to $L$, the
coefficients are $1$-forms and the second equation is obtained by acting by
$J$ on the first equation. Moreover, the metric character of $\nabla^{0}$
implies
$$\lambda_{i}^{j}+\lambda_{j}^{i}=0,\;b^{j}_{i}=b^{i}_{j}.\leqno{(3.14)}$$
The coefficients $\lambda_{i}^{j}$ are the local connection forms of the
metric connection $\nabla^{1}$ induced by $\nabla^{0}$ in $L$, and
$b_{i}^{j}$ are the local components of the {\em second fundamental form}
of $L$.

The connection $\nabla^{1}$ extends to $TM$ by putting
$$\nabla^{1}_{X}JZ=J\nabla^{1}_{X}Z\hspace{1cm}
(Z\in\Gamma TL),$$ and it
becomes the connection $\nabla^{1}$ which is required in (3.5).

In order to use (3.5) we need local unitary bases, and we may use the
bases (3.1) where $e_{i}$ are those of (3.13). Accordingly, we get
$$\nabla^{0}\epsilon_{i}=(\lambda_{i}^{j}+\sqrt{-1}b_{i}^{j})\epsilon_{j},
\;\nabla^{1}\epsilon_{i}=\lambda_{i}^{j}\epsilon_{j}.\leqno{(3.15)}$$
The difference tensor has the components
$$\alpha_{i}^{j}=-\sqrt{-1}b_{i}^{j}.\leqno{(3.16)}$$
The curvature variation can also be obtained by a technical
computation \cite{V1}. In particular, we get
\proclaim 3.1 Proposition. The first Maslov class $\mu_{1}(L)$ is
represented by the differential $1$-form $(1/2\pi)b_{i}^{i}$.\par
This is a generalization of a result due to J. M. Morvan in ${\bf R}^{2n}$
\cite{Mor}, where Proposition 3.1 yields a nice
relationship between the first Maslov
class and the mean curvature vector of a Lagrangian
submanifold $L\subseteq{\bf
R}^{2n}$.

An interesting situation where the calculations above can be used is that
of a symplectic manifold $(M,\omega)$ endowed with a transversally
Riemannian Lagrangian foliation ${\cal V}$. The tangent structure is
provided by Proposition 1.2, and the natural connection $\nabla^{0}$
to be used is that for which $\nabla^{0}/_{V'}$ is the Levi-Civita
connection of the transversal metric of the foliation. In particular, we
can apply these choices to the case of a Lagrangian submanifold of a
tangent bundle $TN$ with a global Lagrangian symplectic structure. What we
will get is a translation of known calculations on cotangent bundles
via a Legendre transformation \cite{{V1},{V6}}.
\section{Locally Lagrangian Poisson Manifolds}
The aim of this section is to suggest an open problem namely, the study of
Poisson manifolds such that their symplectic leaves are l.L.s manifolds.
Besides, we would also like to have the l.L.s. structure of the leaves vary
smoothly, in a reasonable sense. This leads to the following
\proclaim 4.1 Definition. {\rm A {\em locally Lagrangian Poisson (l.L.P.)
manifold} is a triple $(M,P,S)$ where $P$ is a Poisson bivector field
on $M$, and
$S\in\Gamma End\,TM$ and satisfies the properties:\\
$$P(\alpha,\beta\circ S)=P(\beta,\alpha\circ S),\leqno{(4.1)}$$
$$P(\alpha\circ S,\beta\circ S)=0,\leqno{(4.2)}$$
$$rank_{x}S/_{im\,\sharp_{P}}=\frac{1}{2}rank_{x}P,\leqno{(4.3)}$$
$${\cal N}_{S}(X,Y)=0,\hspace{5mm}\forall X,Y\in\Gamma(im\,\sharp_{P}),
\leqno{(4.4)}$$
where $\alpha,\beta\in\Gamma T^{*}M$, $x\in M$,
$\sharp_{P}:T^{*}M\rightarrow TM$ is defined by
$<\sharp_{P}\alpha,\beta>=P(\alpha,\beta)$, and ${\cal N}_{S}$ is the
Nijenhuis tensor (1.7). If all these conditions, with the exception of
(4.4) hold, $(M,P,S)$ is an {\em almost l.L.P. manifold}.}\par
The name is justified by
\proclaim 4.2 Proposition. The symplectic leaves of a l.L.P. manifold are
locally Lagrangian symplectic manifolds.\par
\noindent{\bf Proof.} Formula (4.1) also reads
$$S\sharp_{P}\alpha=-\sharp_{P}(\alpha\circ S),\leqno{(4.5)}$$
and this shows that the tangent spaces of the symplectic leaves are
$S$-invariant. From (4.1), (4.2), it follows that for
$F:=S/_{im\,\sharp_{P}}$, $F^{2}=0$ which, together with (4.3) and (4.4),
shows that $F$ defines a tangent structure on every symplectic leaf ${\cal
S}$ of $P$. Furthermore, if we look at the symplectic structure
$$\omega_{\cal S}(\sharp _{P}\alpha,\sharp_{P}\beta):=-P(\alpha,\beta)
\leqno{(4.6)}$$
of the leaf ${\cal S}$, we see that (4.1) implies the compatibility
condition (1.10). Q.e.d.

An easy example is provided by a manifold
$$M=TN^{n}\times_{f}TN^{n},\leqno{(4.7)}$$
which is the fibered product of two copies of a tangent bundle. In this
case, $M$ has an atlas of local coordinates $(x^{i},y^{i},z^{i})_{i=1}^{n}$
with the coordinate transformations
$$\tilde x^{i}=\tilde x^{i}(x^{j}),\;\tilde y^{i}=
\frac{\partial \tilde x^{i}}{\partial x^{j}}y^{j},\;\tilde z^{i}=
\frac{\partial \tilde x^{i}}{\partial x^{j}}z^{j}, \leqno{(4.8)}$$
and we get a tensor field $S\in\Gamma End\,TM$ if we ask
$$S(X_{i})=0, S\left(\frac{\partial}{\partial y^{i}}\right)=
\frac{\partial}{\partial z^{i}},\;S\left(\frac{\partial}{\partial z^{i}}
\right)=0,
\leqno{(4.9)}$$ where
$$X_{i}=\frac{\partial}{\partial x^{i}}+t^{j}_{i}
\frac{\partial}{\partial z^{j}}\leqno{(4.10)}$$
are chosen such as to span a fixed transversal distribution of the fibers
of the second copy of $TN$ in (4.7).

The definition of $S$ is inspired by (1.3), and it is easy to check that
$S$ is invariant by (4.8) and satisfies
$$S^{2}=0,\;{\cal N}_{S}=0.\leqno{(4.11)}$$

Now, let $P$ be a $2$-contravariant non degenerate symmetric tensor field
on $N$. Then,
$$\Pi=P^{ij}\frac{\partial}{\partial y^{i}}\wedge
\frac{\partial}{\partial z^{j}}\leqno{(4.12)}$$
is a Poisson bivector field on $M$ with the symplectic leaves given by the
fibers of the natural fibration $M\rightarrow N$, and the conditions
(4.1)-(4.4) are satisfied.

Therefore, $(M,\Pi,S)$ is a l.L.P. manifold. Moreover, if we define the
covariant tensor $\Lambda$ such that $\Lambda_{ij}P^{jk}=\delta_{i}^{k}$,
the symplectic structure of the leaf over
any fixed point $x\in M$ is $\Lambda_{ij}(x)
dz^{i}\wedge dy^{j}$, and the leaf has the global Lagrangian function
${\cal L}=\Lambda_{ij}(x)z^{i}y^{j}$, which also is a global function on
$M$.

In agreement with this example we give
\proclaim 4.3 Definition. {\rm A l.L.P.manifold $(M,P,S)$ is a {\em
globally Lagrangian Poisson (g.L.P.) manifold} if there exists a global
function ${\cal L}\in C^{\infty}(M)$ such that its restriction to the
symplectic leaves of $P$ is a global Lagrangian function
of the induced symplectic structure.}\par
\proclaim 4.4 Remark. If $N$ is a locally affine manifold with the affine
local coordinates $(x^{i})$, the manifold (4.7) can also be seen as
$M=T^{(2)}N$, the second order osculating
bundle of $N$ (i.e., the bundle of the second order jets at $0$ of the
mappings in $C^{\infty}({\bf R},M)$).\par
Now, the following question is natural: on a tangent bundle $TN$, find all
the Poisson structures $P$ such that $(M,P,S)$ where $S$ is the canonical
tangent structure of $TN$ is a l.L.P. or a g.L.P. manifold.

If we use the coordinates of (1.3), it follows
$$dq^{i}\circ S=0,\;du^{i}\circ S=dq^{i},\leqno{(4.13)}$$
and (4.2) with $\alpha=du^{i},\beta=du^{j}$ yields $P(df,dg)=0$ for all
$f,g\in C^{\infty}(M)$. This means that $P$ must be a zero-related
structure i.e., such that $\pi:(TN,P)\rightarrow(N,0)$ is a Poisson mapping
\cite{V2}. Furthermore, again using (4.13), we see that (4.1) reduces to
$$P(dq^{i},du^{j})=P(dq^{j},du^{i}).$$
Hence, $P$ must be of the form
$$P=P^{ij}\frac{\partial}{\partial q^{i}}\wedge\frac{\partial}{\partial
u^{j}}+\frac{1}{2}A^{ij}\frac{\partial}{\partial u^{i}}
\wedge\frac{\partial}{\partial u^{j}},\leqno{(4.14)}$$
where $P^{ij}=P^{ji}$, $A^{ij}=-A^{ji}$.

It is an easy consequence of (4.13), (4.14) that $$S(im\,\sharp_{P})=
span\left\{P^{ij}\frac{\partial}{\partial u^{j}}\right\}$$ hence,
$rank\,S=rank(P^{ij})$. Therefore, the problem becomes that of finding the
Poisson bivectors fields which satisfy the condition
$$rank\,P=2rank(P^{ij})\leqno{(4.15)}$$
at each point of $TN$.
 \vspace*{1cm}
{\small Department of Mathematics, \\}
{\small University of Haifa, Israel.\\}
{\small E-mail: vaisman@math.haifa.ac.il}
\end{document}